\documentclass[a4paper,10pt]{scrartcl}

\usepackage[left=3.5cm, right=3.5cm, top=2.3cm, bottom=2.1cm]{geometry}
\usepackage[utf8]{inputenc}
\usepackage{babel}
\usepackage[T1]{fontenc}
\usepackage{amsmath}
\usepackage{amssymb}
\usepackage{multicol}
\usepackage{graphicx}
\usepackage{color}
\usepackage{lmodern}
\usepackage{float}
\usepackage{caption}
\DeclareMathSizes{10}{10}{8.3}{7}

\title{Cut covers of acyclic digraphs}
\author{Maximilian Krone}
\date{October 9, 2024}

\newcounter{ThNr}[subsection]
\DeclareRobustCommand{\newTh}[1]{\refstepcounter{ThNr}\arabic{ThNr} \label{#1}}

\begin{document}
\parskip2.5explus1exminus1ex
\parindent0mm

\begin{center}
%\textbf{\huge{The cut cover problem \\in acyclic digraphs}}\\
\textbf{\huge{Cut covers of acyclic digraphs}}\\
\vspace*{5mm}
{\large Maximilian Krone}

Technische Universität Ilmenau

October 9, 2024 \vspace*{3mm}
\end{center}

\begin{abstract}
A \textit{cut} in a digraph $D=(V,A)$ is a set of arcs $\{uv \in A: u\in U, v\notin U\}$, for some $U\subseteq V$. It is known that the arc set $A$ is covered by $k$ cuts if and only if it admits a $k$-coloring such that no two consecutive arcs $uv, vw$ receive the same color. Alon, Bollobás, Gyárfás, Lehel and Scott (2007) observed that every acyclic digraph of maximum indegree at most $\binom{k}{\lfloor k/2 \rfloor}-1$ is covered by $k$ cuts. 

We prove that this degree condition is best possible (if an enormous outdegree is allowed). Notably, for $k\geq 5$, powers of directed paths do not suffice as extremal examples. Instead, we locate the maximum $d$ such that the $d$-th power of an arbitrarily long directed path is covered by $k$ cuts between $(1-o(1)) \frac{1}{e} 2^k$ and $\frac{1}{2}2^k-2$.

Let $k\geq 3$ and $D$ be an acyclic digraph that is not covered by $k$ cuts. We prove that the decision problem whether a digraph that admits a homomorphism to $D$ is covered by $k$ cuts is NP-complete. If $k=3$ and $D$ is the third power of the directed path on 12 vertices, then even the restriction to planar digraphs of maximum indegree and outdegree $3$ holds.
\end{abstract}

\hrulefill
\section*{Introduction}

Let $D=(V,A)$ be a digraph and $k\in \mathbb{N}$. The following properties are equivalent \cite{Maxcuts, arc colorings}:
    \begin{enumerate}
        \item[(i)] $A$ is covered by $k$ \textit{cuts} $\{uv \in A : u\in U_a, v\notin U_a\}$, $U_1,\dots,U_k \subseteq V$. (We say, $D$ has a $k$-\textit{cut cover}.) \vspace*{-2mm}
        \item[(ii)] There is a $k$-coloring $c: A \to \{1,\dots,k\}$ such that no two consecutive arcs $uv, vw$ receive the same color. \vspace*{-2mm}
        \item[(iii)] There is a labeling $C: V \to \mathcal{P}\left( \{1,\dots,k\} \right)$ (the \textit{characteristic sets}) such that, for every arc $uv \in A$, there is an $a$ with $a \in C(u), a \notin C(v)$, that is, $C(u)\not \subseteq C(v)$.
    \end{enumerate}
    
\begin{figure}[H]
\begin{center}
	\includegraphics[height=18mm]{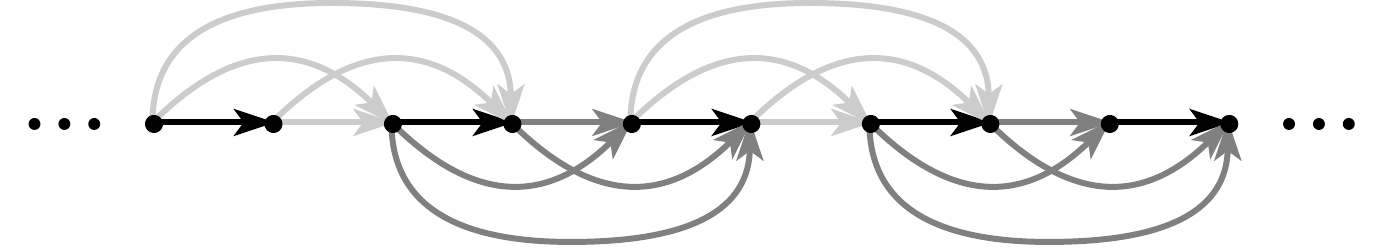}
	\vspace*{-2mm} \caption{A 3-coloring of the arc set as in (ii). In this example, each color sets is a cut itself, while in general, it is only contained in some cut. The characteristic sets can be chosen as $\dots\, ,\,\{1,3\}\, ,\,\{1\}\, ,\,\{2,3\}\, ,\,\{2\}\, ,\,\{1,3\}\, ,\,\{1\}\, ,\,\{2,3\}\, ,\,\{2\}\, ,\,\{1,3\}\, ,\,\{1\}\, ,\,\dots$}
	\label{example}
\end{center}
\end{figure} \vspace{-7mm}

The properties (i) and (ii) are easier to motivate while property (iii) is most useful for the analysis. Hence, we will always talk about cut covers in the theorems but use characteristic sets in their proofs. Clearly, the characteristic sets of a $k$-cut cover define a $2^k$-coloring of the vertex set \cite{Maxcuts, arc colorings}.

The following inverse observation is even more important. Every digraph that admits a $\binom{k}{\lfloor k/2 \rfloor}$-coloring is covered by $k$ cuts by choosing the characteristic set of every vertex depending on its color: As the colors, use the subsets of $\{1,\dots,k\}$ of the size $\lfloor k/2 \rfloor$, which are pairwise incomparable, so this yields the characteristic sets of a $k$-cut cover.

Every acyclic digraph of maximum indegree at most $\binom{k}{\lfloor k/2 \rfloor}-1$ has a greedy $\binom{k}{\lfloor k/2 \rfloor}$-coloring by considering the vertices in an acyclic ordering. Hence, it is covered by $k$ cuts, which was observed by Alon, Bollobás, Gyárfás, Lehel and Scott in \cite{Maxcuts}. 

In Section \ref{indegree}, we prove that this degree condition is best possible by using a construction that demands vertices of enormous outdegree. One could hope that also powers of directed paths suffice as extremal examples. As it turns out, they only do for $k\leq 4$. More precisely, the maximum $d$ such that the $d$-th power of an arbitrarily long directed path is covered by $k$ cuts is located between $(1-o(1)) \tfrac{1}{e} 2^k$ and $\tfrac{1}{2}2^k-2$, which can be derived from \cite{self}.

In Section \ref{or}, we generalize the degree condition: For any integers $\Delta^-, \Delta^+ \geq -1$, let $\mathcal{D}(\Delta^-,\Delta^+)$ be the class of all acyclic digraphs in which every vertex has indegree at most $\Delta^-$ \textit{or} outdegree at most $\Delta^+$. Let $c(\Delta^-,\Delta^+)$ be the smallest $k$ such that every digraph in $\mathcal{D}(\Delta^-,\Delta^+)$ is covered by $k$ cuts. We are able to limit $c(\Delta^-,\Delta^+)$ to at most two explicit values for each pair $(\Delta^-,\Delta^+)$. For small instances, we also determine the exact values.

Finally, we examine the algorithmic complexity of the decision problem whether a digraph is covered by $k$ cuts. For $k\leq 2$, the decision can be done efficiently. For $k\geq 3$, it is known that the problem is NP-complete for the class of symmetric digraphs, since a symmetric digraph is covered by $k$ cuts if and only if it is $\binom{k}{\lfloor k/2 \rfloor}$-colorable (see Appendix) \cite{symmetric}. 

In Section \ref{Complexity}, we prove that the decision problem is also NP-complete for the class of acyclic digraphs, even under the condition that the input digraph admits a homomorphism to some fixed acyclic digraph that is not covered by $k$ cuts. On the other hand, a homomorphism to a digraph that is covered by $k$ cuts is sufficient for a $k$-cut cover. Powers of directed paths are of particular interest since the existence of a homomorphism to these can be tested efficiently. The decision problem whether an acyclic digraph is covered by $3$ cuts is NP-complete, even for planar digraphs of maximum indegree and outdegree $3$ that admit a homomorphism to the third power of the directed path on 12 vertices.

For shorter notations, we define the central binomial coefficient
$$M(k) := \tbinom{k}{\lfloor k/2 \rfloor} = \tfrac{k!}{\lfloor k/2 \rfloor ! \lceil k/2 \rceil !} = \tbinom{k}{\lceil k/2 \rceil}.$$ 
By Sperner's theorem, $M(k)$ is the size of the largest family of pairwise incomparable subsets of $\{1,\dots,k\}$ \cite{Sperner}. The first values are $\left(M(k)\right)_{k\geq 1} = (1,2,3,6,10,20,35,70,\dots)$. 
The bounds $\frac{\sqrt{2}}{2} \frac{2^k}{\sqrt{k+1}} \ \leq\  M(k)\ \leq\ \frac{\sqrt{3}}{2} \frac{2^k}{\sqrt{k+1}}\ $ can be proved simply by induction (see Appendix). %Asymptotically,\ $M(k)\sim \sqrt{\frac{2}{\pi}}\ \frac{2^k}{\sqrt{k}}$, by Stirling's approximation.

\section{Cut covers of acyclic digraphs of bounded indegree} \label{indegree}

In \cite{Maxcuts}, Alon, Bollobás, Gyárfás, Lehel and Scott observed that every acyclic digraph with maximum outdegree at most $M(k)-1$ is covered by $k$ cuts, since it is (greedily) $M(k)$-colorable. Clearly, we have the same condition for the maximum indegree. We are going to show that the bound $M(k)-1$ is best possible. But to bound the indegree, the construction uses vertices of enormous outdegree.

The complete digraph on $M(k)+1$ vertices shows that not every (non-acyclic) digraph with both maximum indegree and outdegree $M(k)$ is covered by $k$ cuts, since each two characteristic sets have to be incomparable. Unfortunately, it remains open whether this statement is also true for acyclic digraphs. At least, we can transfer it for $k\leq 4$, by considering powers of directed paths as the simplest class of acyclic digraphs with both bounded indegree and outdegree.

The $d$-th power $P_n^d$ of the directed path on $n$ vertices is the digraph on vertex set $\{1,\dots,n\}$ (or $\{0,\dots,n-1\}$) that contains an arc $ij$ if and only if $j-i \in \{1,\dots,d\}$.
For $k\geq 1$, let $\delta(k)$ be the maximum $d$ such that, for all $n\in \mathbb{N}$, $P_n^d$ is covered by $k$ cuts.

Clearly, $\delta(k)\leq 2^k-1$, since the transitive tournament $P_{2^k+1}^{2^k}$ has no $2^k$-coloring and hence no $k$-cut cover. On the other hand, since $P_n^{M(k)-1}$ is $M(k)$-colorable, $\delta(k) \geq M(k)-1$. \newpage

The definition of $\delta(k)$ clearly handles the same object as \cite{self}. We derive the following two theorems, as well as the conjecture. 

\textbf{Theorem \ref{indegree}.\newTh{block path}\ -\ A class of digraphs close to powers of directed paths\\}
Let $k\geq 2$. Consider the acyclic digraph $D$ on vertex set $\{1,\dots,l\}\times \{1,\dots,\frac{1}{4}2^k\}$ and arc set $\lbrace (r,i)(r',i') : r=r', i<i' \text{ or } r+1=r'\rbrace$. Let $D'$ be build from $D$ by adding a new vertex $v$ and the arcs $(\lfloor l/2 \rfloor, i)v,\ v(\lfloor l/2 \rfloor+1, i)$, for all $i\in \{1,\dots,\frac{1}{4}2^k\}$. 
\begin{enumerate}
	\item[(i)] For all $l\in \mathbb{N}$, $D$ is covered by $k$ cuts. \vspace*{-2mm}
	\item[(ii)] If $l\geq 2^{(k+3)/2}$, then $D'$ is not covered by $k$ cuts. For every vertex in $D'$, the sum of its indegree and outdegree is at most $\frac{3}{4}2^k$.\\
\end{enumerate}

For $k=3$, $D$ and its cut cover are shown in Figure \ref{example}. The construction can easily be generalized. Since $D$ can contain the $\frac{1}{4}2^k$-th power of an arbitrarily long directed path, $\delta(k) \geq \frac{1}{4}2^k$. On the other hand, since $D'$ is contained in the $\frac{1}{2}2^k$-th power of some directed path, $\delta(k) \leq \frac{1}{2}2^k -1$. Hence, the behavior of $\left(\delta(k)\right)_{k\geq 1}$ is more similar to the simple upper bound $\delta(k)\leq 2^k-1$. With some more ideas, the following results are achieved.

\textbf{Theorem \ref{indegree}.\newTh{delta}\ -\ Powers of directed paths}
\begin{enumerate} \vspace{-2mm}
	\item[(i)] For $k\geq 4$, \ $\delta(k) \geq \tfrac{5}{16} 2^k$. \vspace*{-2mm}
	\item[(ii)] If $r$ is a power of $2$ and $k\geq r \log_2 r$, then $\delta(k) \geq \left( 1-\tfrac{1}{r} \right)^r 2^k \geq \left( 1-\tfrac{1}{r} \right) \tfrac{1}{e} 2^k$. \vspace*{-3mm}
	\item[(iii)] $P^d_{2^k+1}$ is covered by $k$ cuts if and only if  $d \leq 2^k -2^{\lfloor k/2 \rfloor} - 2^{\lceil k/2 \rceil} +1$. \vspace*{-2mm}
	\item[(iv)] $P^3_n$ is covered by $3$ cuts if and only if $n\leq 10$, in particular $\delta(3)=2$. \vspace*{-1mm}
	\item[(v)] $P^6_n$ is covered by $4$ cuts if and only if $n\leq 24$, in particular $\delta(4)=5$. \vspace*{-1mm}
	\item[(vi)] For $k\geq 3$, \ $\delta(k) \leq \frac{1}{2}2^k -2$.\\
\end{enumerate}

The author conjectures in \cite{self} that (ii) is essentially best possible.

\textbf{Conjecture\ -\ \boldmath $\delta(k)=(1-o(1))\tfrac{1}{e}2^k$.} \vspace*{-2mm}
$$\delta(k) \leq \max\left\{d\in \mathbb{N}: 1\leq \frac{d^d}{(d+1)^{d+1}}2^k\right\} 
< \frac{1}{e}2^k\ .$$
By the statements (iii)-(v), the trivial lower bound $\delta(k) \geq M(k)-1$ is attained for $k\leq 4$. But by (i), it is not for $k\geq 5$. Therefore, to prove that, for all $k\in\mathbb{N}$, not every acyclic digraph of maximum indegree $M(k)$ is covered by $k$ cuts, we need a cleverer construction. Unfortunately, we also have to be a lot coarser and use vertices of outdegree $M(k)^{\omega(1)M(k)}$.

\textbf{Theorem \ref{indegree}.\newTh{Th indegree}\ -\ \boldmath Not every acyclic digraph of maximum indegree $M(k)$ is covered by $k$ cuts \\ }
For $d,k\in \mathbb{N}$, there exists an acyclic digraph $D=(V,A)$ with the following properties.
\begin{enumerate}
\item[(i)] $D$ has maximum indegree $d$. \vspace*{-2mm}
\item[(ii)] Every vertex and its in-neighborhood induce a transitive tournament. \vspace*{-2mm}
\item[(iii)] $D$ admits a homomorphism to the acyclic tournament on $2^k+1$ vertices. \vspace*{-2mm}
\item[(iv)] If $D$ is covered by $k$ cuts, then there is some vertex of indegree $d$ for which the characteristic sets of itself and its in-neighbors are pairwise incomparable.
\end{enumerate}
In particular, if $d \geq M(k)$, by Sperner's theorem, there is no family of $d+1$ pairwise incomparable subsets of $\{1,\dots,k\}$ as in (iv). Hence, $D$ is not covered by $k$ cuts.

\textbf{Proof. \ } We give an explicit construction of $D=(V,A)$. For any vertex $v \in V$, let \linebreak $N^-(v) = \{u \in V: uv \in A\}$ be its in-neighborhood. 

We construct $V$ iteratively in layers $V_0,\dots,V_{2^k}$. For $i=0,\dots,d$, let $V_i = \{v_i\}$ with $N^-(v_i)=\{v_0,\dots,v_{i-1}\}$.
For $i=d+1,\dots,2^k$, let $V_i = \{ w_{vu}:\ v\in V_{i-1},\ u \in N^-(v) \}$ with $N^-(w_{vu}) = \{v\} \cup N^-(v)\setminus \{u\}$.

Properties (i) and (ii) are true for $v_0,\dots,v_d$. Inductively, all vertices from layers with higher index have indegree $d$ and they form transitive tournaments with their in-neighborhood. So (i) and (ii) hold.
$D$ is acyclic and (iii) holds since all arcs go from layers with lower to layers with higher index.

The last open claim is (iv). Fix a $k$-cut cover of $D$. Assume that there is no vertex $v$ of indegree $d$ for which the characteristic sets of $\{v\} \cup N^-(v)$ are pairwise incomparable.

For any vertex set $U \subseteq V$, we define the family of sets 
$$\mathcal{B}(U) := \{B\subseteq \{1,\dots,k\}: B\supseteq C(u) \text{ for some } u\in U\}.$$
Let $v\in V$. The main observation is that $C(v) \notin \mathcal{B}(N^-(v))$, but clearly $C(v) \in \mathcal{B}(\{v\} \cup N^-(v))$. Hence, $\mathcal{B}(N^-(v)) \subsetneq \mathcal{B}(\{v\} \cup N^-(v))$. 

Iteratively and depending on the $k$-cut cover, we choose a sequence $v_0,\dots,v_{2^k}$ of vertices such that $\mathcal{B}(N^-(v_0)) \subsetneq \mathcal{B}(N^-(v_1)) \subsetneq \dots \subsetneq \mathcal{B}(N^-(v_{2^k}))$. This implies $|\mathcal{B}(N^-(v_{2^k}))| \geq 2^k$, which contradicts $C(v_{2^k}) \notin \mathcal{B}(N^-(v_{2^k}))$.
For $i\leq d$, let $v_i$ be the unique vertex in $V_i$, so indeed $\mathcal{B}(N^-(v_{i-1})) \subsetneq \mathcal{B}(\{v_{i-1}\} \cup N^-(v_{i-1})) = \mathcal{B}(N^-(v_i))$.

For $i > d$, by the assumption, among the characteristic sets of $\{v_{i-1}\} \cup N^-(v_{i-1})$, there are two comparable sets $C(u)\supseteq C(w)$, $u \neq w$. Since no arc $wu$ is allowed, we have $u\neq v_{i-1}$, so $u\in N^-(v_{i-1})$. Now, choose $v_i = w_{v_{i-1}u} \in V_i$. For every set $B \supseteq C(u)$, it also holds $B \supseteq C(w)$. This ensures that $\mathcal{B}(N^-(v_{i-1})) \subsetneq \mathcal{B}(\{v_{i-1}\} \cup N^-(v_{i-1})) = \mathcal{B}(\{v_{i-1}\} \cup N^-(v_{i-1})\setminus \{u\}) =  \mathcal{B}(N^-(w_{v_{i-1}u})) = \mathcal{B}(N^-(v_i))$.  \hfill $\Box$ \\

One of the original questions remains open: For every $k$, does there exist an acyclic digraph with both indegree and outdegree bounded by $M(k)$ (or $o(1)2^k$) that is not covered by $k$ cuts? The first unsolved case is $k=5$. For this case, even with some slight improvements, the construction from Theorem \ref{indegree}.\ref{Th indegree} demands vertices of outdegree $\approx 10^{19}$.

\section{Cut covers of acyclic digraphs with vertices of bounded indegree and vertices of bounded outdegree} \label{or}

For any integers $\Delta^-, \Delta^+ \geq -1$, let $\mathcal{D}(\Delta^-,\Delta^+)$ be the class of all acyclic digraphs in which every vertex has indegree at most $\Delta^-$ \textit{or} outdegree at most $\Delta^+$. Let $c(\Delta^-,\Delta^+)$ be the smallest $k$ such that every digraph in $\mathcal{D}(\Delta^-,\Delta^+)$ is covered by $k$ cuts (similarly to \cite{Maxcuts}).

Clearly, $c(\Delta^-,\Delta^+)$ is symmetric and monotone in $\Delta^-$ and $\Delta^+$.
The case $\Delta^+ =-1$ generates the class $\mathcal{D}(\Delta^-,-1)$ of all acyclic digraphs of maximum indegree at most $\Delta^-$. Hence, we can derive the values of $c(\Delta^-,-1)$ for all $\Delta^- \in \mathbb{N}$ from Theorem \ref{indegree}.\ref{Th indegree}. 
Degenerate cases are the class $\mathcal{D}(-1,-1)$, which consists of the empty digraph only, and the classes $\mathcal{D}(-1,0)=\mathcal{D}(0,-1)$, which consist of all digraphs with empty arc set. The class $\mathcal{D}(0,0)$ consists of those digraphs whose arc set is a cut itself: All arcs go from vertices with indegree 0 to vertices with outdegree 0. \enlargethispage{12mm}

\textbf{Theorem \ref{or}.\newTh{Three bounds}\ -\ Three bounds on \boldmath $c(\Delta^-,\Delta^+)$} 
\begin{enumerate}
	\item[(i)] If $\Delta^-+\Delta^+ \leq M(k) -2$, then $c(\Delta^-,\Delta^+) \leq k$ \cite{Maxcuts}.
	\item[(ii)] If $\max \{\Delta^-,\Delta^+\} > M(k)-1$, then $c(\Delta^-, \Delta^+) > k$.
	\item[(iii)] If $\Delta^-+\Delta^+ > M(k+1)-2$, then $c(\Delta^-, \Delta^+) > k$.
\end{enumerate}

\begin{figure}[H]
\begin{center}
	\vspace*{-4mm}
	\includegraphics[height= 97mm]{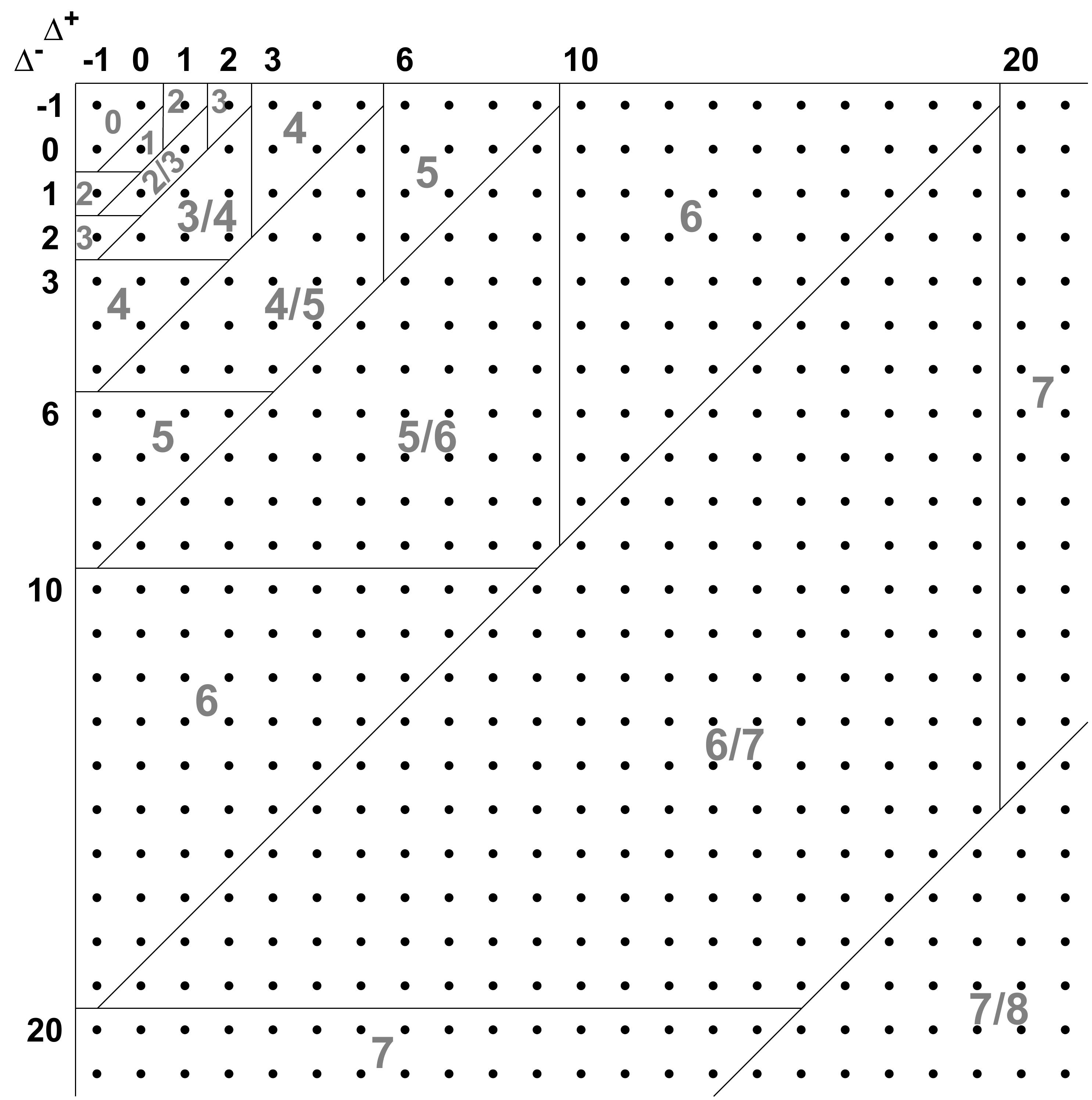}
	\caption{The possible values of $c(\Delta^-,\Delta^+)$ according to Theorem \ref{or}.\ref{Three bounds}. Each point corresponds to a pair $(\Delta^-,\Delta^+)$. Bound (ii) is stronger than bound (iii) in most cases, except for even $k$ and $\Delta^-\approx \Delta^+$.}
	\label{threebounds}
\end{center}
\end{figure} \vspace{-8mm} 

\textbf{Proof. \ } (i) -\ Let $D=(V,A)\in \mathcal{D}(\Delta^-,\Delta^+)$. Let $V^-$ be the set of all vertices with indegree at most $\Delta^-$. Then, all vertices in $V^+ = V \setminus V^-$ have outdegree at most $\Delta^+$. The subdigraphs of $D$ that are induced by $V^-$ and $V^+$ are (greedily) $(\Delta^-+1)$- and $(\Delta^++1)$-colorable, respectively. Note that this is also true for $\Delta^-=-1$ or $\Delta^+=-1$. By using disjoint color sets, this yields a vertex coloring of $D$ with $\Delta^-+\Delta^++2\leq M(k)$ colors. Hence, $D$ is covered by $k$ cuts \cite{Maxcuts}. 

(ii) -\ By Theorem \ref{indegree}.\ref{Th indegree}, there exists some acyclic digraph of maximum indegree or, by inverting all arcs, maximum outdegree $M(k)$ that is not covered by $k$ cuts. This digraph is in $\mathcal{D}(M(k),-1)$ or $\mathcal{D}(-1,M(k))$, respectively.

(iii) -\ Let $k$ be odd. It follows from $\Delta^-+\Delta^+ > M(k+1)-2 = 2M(k)-2$ that $\max\{\Delta^-,\Delta^+\} > M(k)-1$, so by (ii), $c(\Delta^-, \Delta^+) >k$. 

To prove the claim for even $k$, we have to use the construction of Theorem \ref{indegree}.\ref{Th indegree} twice.
Let $D^-=(V^-,A^-)$ be an acyclic digraph with maximum indegree $\Delta^-$ and the property that, for every $k$-cut cover, there appears a family $\mathcal{B}^-$ of $\Delta^-+1$ pairwise incomparable characteristic sets. By inverting all arcs in a similar construction, we obtain an acyclic digraph $D^+=(V^+,A^+)$  with maximum outdegree $\Delta^+$ and the property that, for every $k$-cut cover, there appears a family $\mathcal{B}^+$ of $\Delta^++1$ pairwise incomparable characteristic sets. 

We assume that $V^-$ and $V^+$ are disjoint, and we unite both digraphs to obtain an acyclic digraph $D=(V,A)$ with $V = V^- \cup V^+$ and $A = A^- \cup A^+ \cup (V^- \times V^+)$. 
In $D$, all vertices in $V^-$ still have indegree at most $\Delta^-$ and all vertices in $V^+$ still have outdegree at most $\Delta^+$. That is, $D \in \mathcal{D}(\Delta^-,\Delta^+)$.

Assume that $D$ is covered by $k$ cuts and consider its characteristic sets. Since $V^- \times V^+ \subseteq  A$, the two families $\mathcal{B}^-$ and $\mathcal{B}^+$ of $\Delta^-+1$ and $\Delta^++1$ pairwise incomparable sets, respectively, are disjoint.
For $i \in \{0,\dots,k\}$, let $b^-_i$ and $b^+_i$ be the number of all $i$-element sets in $\mathcal{B}^-$ and $\mathcal{B}^+$, respectively. 
The inequality of Lubell, Yamamoto, Meshalkin and Bollobás states \vspace*{-1mm}
$$\sum _{i=0} ^k \frac{b^-_i}{\binom{k}{i}} \leq 1 \quad \text{and}\quad \sum _{i=0} ^k \frac{b^+_i}{\binom{k}{i}} \leq 1\ .$$
Let $k=2l$. Since $\mathcal{B}^-$ and $\mathcal{B}^+$ are disjoint, we can bound $b^-_l+b^+_l \leq \binom{2l}{l}$. Using that ${\binom{2l}{i}} \leq {\binom{2l}{l-1}}$ for all $i\neq l$, it follows that
$$ \begin{aligned} 2 \ & \geq \ \sum _{i=0} ^{2l} \frac{b^-_i + b^+_i}{\binom{2l}{i}}\ \geq \ \frac{1}{\binom{2l}{l-1}} \sum_{i=0}^{2l}(b^-_i +b^+_i) -  \left( \frac{1}{\binom{2l}{l-1}} - \frac{1}{\binom{2l}{l}} \right) (b^-_l+b^+_l) \\
 \ & \geq \ \frac{\Delta^-+\Delta^++2}{\binom{2l}{l-1}} - \left( \frac{1}{\binom{2l}{l-1}} - \frac{1}{\binom{2l}{l}} \right) \binom{2l}{l} \
  = \ \frac{\Delta^-+\Delta^++2 - \binom{2l}{l}}{\binom{2l}{l-1}} + 1\ . \end{aligned} $$
Rearranging yields \ $\Delta^-+\Delta^++2 \ \leq \ \binom{2l}{l} + \binom{2l}{l-1} = \binom{2l+1}{l} = M(2l+1) = M(k+1)$. \hfill $\Box$ \\

We now extend the ideas of Theorem \ref{or}.\ref{Three bounds} to identify all pairs $(\Delta^-,\Delta^+)$ with $c(\Delta^-,\Delta^+) \leq k$, for the values $k\leq 4$. For $k=5$, this becomes much more complicated, and we will give none of the proofs here. The pairs $(\Delta^-,\Delta^+)$ with $c(\Delta^-,\Delta^+) \leq 2$ can be identified by using only the following theorem.

\textbf{Theorem \ref{or}.\newTh{2-cut cover}\ -\ 2-cut covers}\\
A digraph is covered by 2 cuts if and only if its subdigraph induced by all vertices with non-zero indegree and outdegree is 2-colorable (which can be tested efficiently).

\textbf{Proof. \ } Let $D=(V,A)$ be a digraph and $W := \{v\in V: d^+(v) \neq 0,\ d^-(v) \neq 0\}$. 
Assume first that there exists a 2-coloring of the subdigraph induced by $W$, $f: W \rightarrow \{1,2\}$ such that $f(u) \neq f(v)$ for every arc $uv \in A$. We choose the characteristic set of any vertex $v$ as
$$
C(v) =  \left\{\begin{aligned} \{1,2\}\ , &\quad d^-(v)=0 \\ \{f(v)\}\ , &\quad d^+(v) \neq 0,\ d^-(v) \neq 0 \\ \emptyset\ , &\quad d^+(v) = 0 \\  \end{aligned}\right.
$$
If both $d^+(v)=d^-(v)=0$, then $C(v)$ can be chosen arbitrarily. Let $uv \in A$. If both $u, v \in W$, we have $C(u) = \{f(u)\} \not\subseteq \{f(v)\} = C(v)$. If $u\notin W$, by the definition of $W$ and $d^+(u)\neq 0$, we get $d^-(u)= 0 \allowbreak \neq d^-(v)$, and so $C(u) = \{1,2\} \not\subseteq C(v)$. Similarly, if $v\notin W$, we get $d^+(u) \neq 0 = d^+(v)$, and so $C(u) \not\subseteq \emptyset = C(v)$. Hence, the characteristic sets indeed describe a $2$-cut cover of $D$.

Assume now that $D$ is covered by $2$ cuts and let $C(v)$, $v \in V$, be the characteristic sets. Let $v \in W$. Then, there exist vertices $u$ and $w$ such that $uv, vw \in A$, and hence $C(u) \not\subseteq C(v) \not\subseteq C(w)$. Therefore, $\{1,2\}\neq C(v)\neq\emptyset$, and we can choose the color $f(v)$ as the unique element of $C(v)$.
Now, for every arc $uv \in A$ with $u,v \in W$, it holds that $\{f(u)\} = C(u) \not\subseteq C(v) = \{f(v)\}$. Hence $f(u) \neq f(v)$, so $f$ is indeed a coloring of the subdigraph induced by $W$. \hfill $\Box$ \\

\begin{figure}[H]
\begin{center}
	\vspace*{-3mm}
	\includegraphics[height= 65mm]{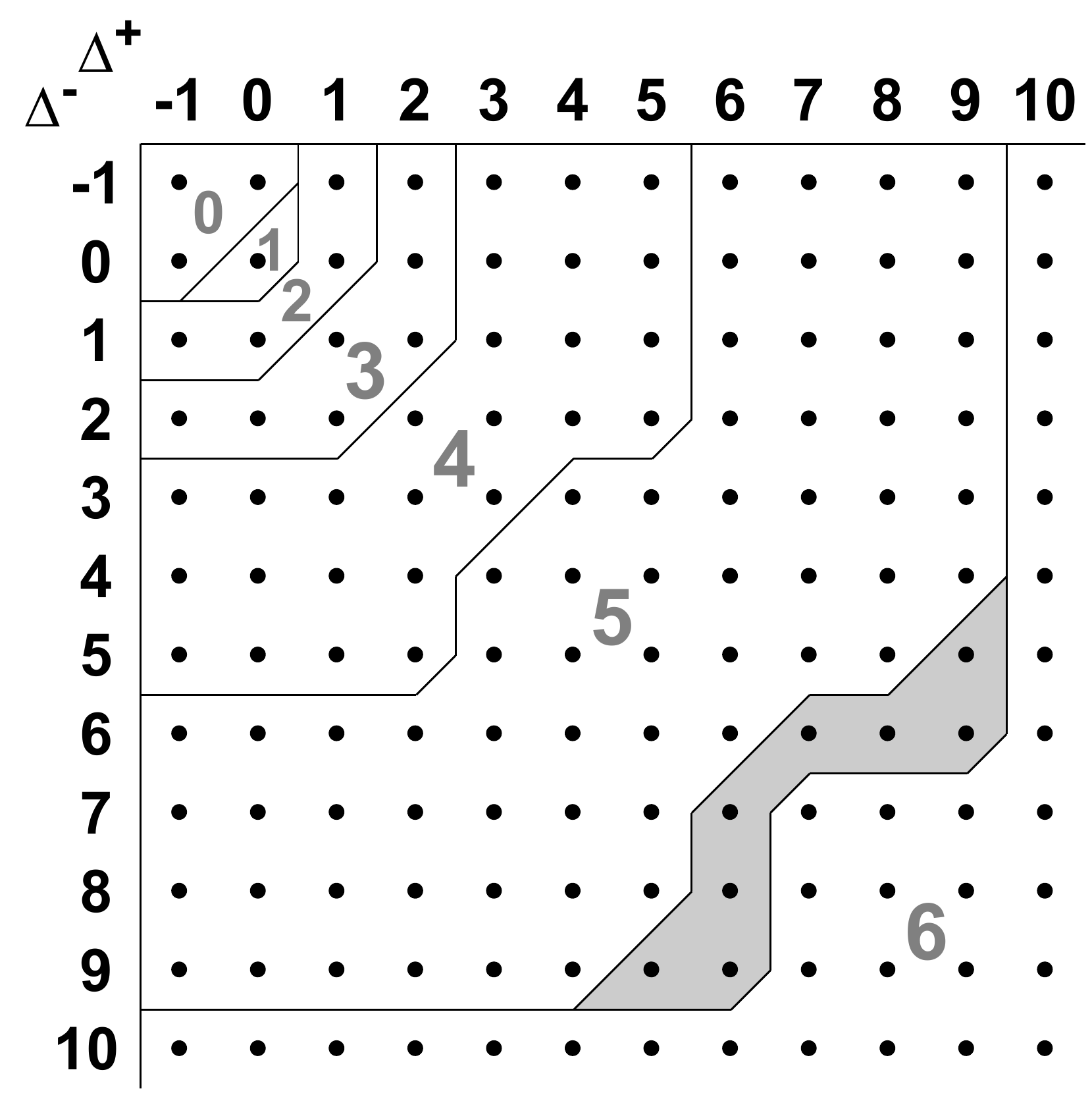} \vspace*{-2mm}
	\caption{The values of $c(\Delta^-,\Delta^+)$. Each point corresponds to a pair $(\Delta^-,\Delta^+)$. The values in the grey area have not been determined yet, but are either $5$ or $6$.}
	\label{values}
\end{center}
\end{figure} \vspace{-6mm} 

We have $c(1,0) \leq 2$, since for every digraph in $\mathcal{D}(1,0)$, the subdigraph induced by the vertices of indegree at most 1 is $2$-colorable and all other vertices have outdegree $0$. 
Both $c(2,-1)>2$ and $c(1,1)>2$ are shown by the digraph with vertices $0,1,2,3,4$ and arcs $01,12,13,23,34$, which is not covered by $2$ cuts, since the vertices $1,2,3$ of non-zero indegree and outdegree form a triangle \cite{Maxcuts}.

\textbf{Theorem \ref{or}.\newTh{characterization3}\ -\ \boldmath Characterization of all pairs $(\Delta^-,\Delta^+)$ with $c(\Delta^-,\Delta^+) \leq 3$}
$$c(2,1)\leq 3\ , \quad c(3,-1)>3\ , \quad c(2,2)>3\ .$$
{\boldmath $c(2,1)\leq 3: \ $\ }  Let $D=(V,A)\in \mathcal{D}(2,1)$, and let $V^-$ be the set of all vertices with indegree at most $2$. Then, all vertices in $V^+=V \setminus V^-$ have outdegree at most $1$.

First, we color $V^-$ greedily with the three colors $\{1,2\}, \{1,3\}, \{2,3\}$, which we use as characteristic sets. This ensures $C(u)\not\subseteq C(v)$ for all $uv\in A$ with $u,v\in V^-$.

Afterwards, we greedily choose the characteristic sets of the vertices in $V^+$ as one of the sets $\{1\},\{2\},\{3\}$. This ensures $C(u)\not\subseteq C(v)$ for all $uv\in A$ with $u\in V^-,v\in V^+$. We proceed as follows.
If the single out-neighbor $v$ of a vertex $u\in V^+$ is also in $V^+$, choose $C(u)\neq C(v)$ arbitrarily. If $v \in V^-$, choose $C(u) = \{1,2,3\} \setminus C(v)$. This ensures $C(u)\not\subseteq C(v)$ for all $uv\in A$ with $u\in V^+$.

{\boldmath $c(3,-1)>3: \ $\ } By Theorem \ref{indegree}.\ref{delta}(iv), $P^3_{11}$ is not covered by $3$ cuts. This is a much smaller example than the digraph that is obtained from Theorem \ref{indegree}.\ref{Th indegree}.

{\boldmath $c(2,2)>3: \ $\ } We are going to prove in Theorem \ref{Complexity}.\ref{NP3}(ii) that the decision problem whether a digraph in $\mathcal{D}(2,2)$ is covered by $3$ cuts is NP-complete. One can derive a simple construction of a digraph in $\mathcal{D}(2,2)$ that is not covered by $3$ cuts from the digraph $P$ that is given in the proof. \hfill $\Box$\\

\textbf{Theorem \ref{or}.\newTh{characterization4}\ -\ \boldmath Characterization of all pairs $(\Delta^-,\Delta^+)$ with $c(\Delta^-,\Delta^+) \leq 4$}
$$c(5,2)\leq 4\ , \quad c(3,3)\leq 4\ , \quad c(6,-1)>4\ , \quad c(4,3)>4\ .$$
{\boldmath $c(5,2)\leq 4: \ $\ } Let $D=(V,A)\in \mathcal{D}(5,2)$, and let $V^-$ be the set of all vertices with indegree at most $5$. Then, all vertices in $V^+=V \setminus V^-$ have outdegree at most $2$.

First, we color $V^-$ greedily with the six colors $B \subseteq \{1,2,3,4\}$ of size $2$, which we use as characteristic sets. This ensures $C(u)\not\subseteq C(v)$ for all $uv\in A$ with $u,v\in V^-$. We add an additional rule to the greedy algorithm: Whenever a vertex $v \in V^+$ has two out-neighbors in $V^-$, and the first one is colored with a set $B$, assign $\{1,2,3,4\}\setminus B$ as a temporary color to $v$, and consider this color in the greedy coloring of $V^-$, that is, when the second out-neighbor of $v$ is colored, it obtains a color different from $\{1,2,3,4\}\setminus B$.

Afterwards, we delete the temporary colors of the vertices in $V^+$ and greedily choose the characteristic sets of the vertices in $V^+$ as one of the sets $\{1\},\{2\},\{3\},\{4\}$. This ensures $C(u)\not\subseteq C(v)$ for all $uv\in A$ with $u\in V^-,v\in V^+$. We proceed as follows.

The characteristic sets of the out-neighbors of $u$ that are in $V^-$ contain only one element. If $u$ has two out-neighbors in $V^+$, we have chosen the second one not disjoint from the first one.
In all cases, there is some $a\in\{1,2,3,4\}$ that is contained in neither of the characteristic sets of the out-neighbors, and we choose $C(u)=\{a\}$. This ensures $C(u)\not\subseteq C(v)$ for all $uv\in A$ with $u\in V^+$.

{\boldmath $c(3,3)\leq 4:\ $\ } Let $D = (V,A) \in \mathcal{D}(3,3)$, let $V^-$ be the set of vertices with indegree at most $3$, and $V^+ = V \setminus V^-$. We greedily color $V^-$ with the four colors $\{1,4\}, \{2,4\}, \{3,4\}, \{1,2,3\}$ such that only vertices with 3 in-neighbors in $V^-$ get color $\{1,2,3\}$. Similarly, we greedily color $V^+$ with the colors $\{1,2\}, \{2,3\}, \{1,3\}, \{4\}$ such that only vertices with 3 out-neighbors in $V^+$ get color $\{4\}$. 

Let $uv \in A$. We have to check $C(u) \not\subseteq C(v)$. If both $u,v \in V^-$ or both $u,v\in V^+$, this is true since both color families that we used for the coloring consist of pairwise incomparable sets. If $u\in V^-$ and $v\in V^+$, this is true since $C(u)\neq C(v)$ and $|C(u)| \geq |C(v)|$. If $u \in V^+$ and $v \in V^-$, the additional property of the coloring yields $|C(u)| = |C(v)| = 2$, so the claim follows from $C(u)\neq C(v)$.

{\boldmath $c(6,-1)>4:\ $\ } By Theorem \ref{indegree}.\ref{delta}(v), $P^6_{25}$ is not covered by $4$ cuts. This is a much smaller example than the digraph that is obtained from Theorem \ref{indegree}.\ref{Th indegree}.

{\boldmath $c(4,3)>4:\ $\ } Let $D^-=(V^-,A^-)$ and $D^+=(V^+,A^+)$ be constructed as in Theorem \ref{or}.\ref{Three bounds}(iii), for $\Delta^-=4$ and $\Delta^+=3$. Consider the following digraph $D\in\mathcal{D}(4,3)$. We build the vertex set from $V^-$, $V^+$ and four more vertices $x_1,x_2,x_3,x_4$. Moreover, we add a set of vertices $W = \{ w_{\{x,y,z\}}: x\in \{x_1,x_2,x_3\}, y,z\in V^+ \}$. 

We build the arc set from $A^-\ \cup\  A^+\ \cup\ \left(V^-\times (W\cup\{x_1,x_2,x_3\}\cup V^+)\right)\ \cup\ (\{x_4\} \times V^+)\ \cup\ \{x_ix_j: 1\leq i < j \leq 4\}$. Moreover, from each $w_{\{x,y,z\}} \in W$, we add arcs to $x,y,z$.

\begin{figure}[H]
\begin{center}
	\vspace*{-3mm}
	\includegraphics[width= 8.5cm]{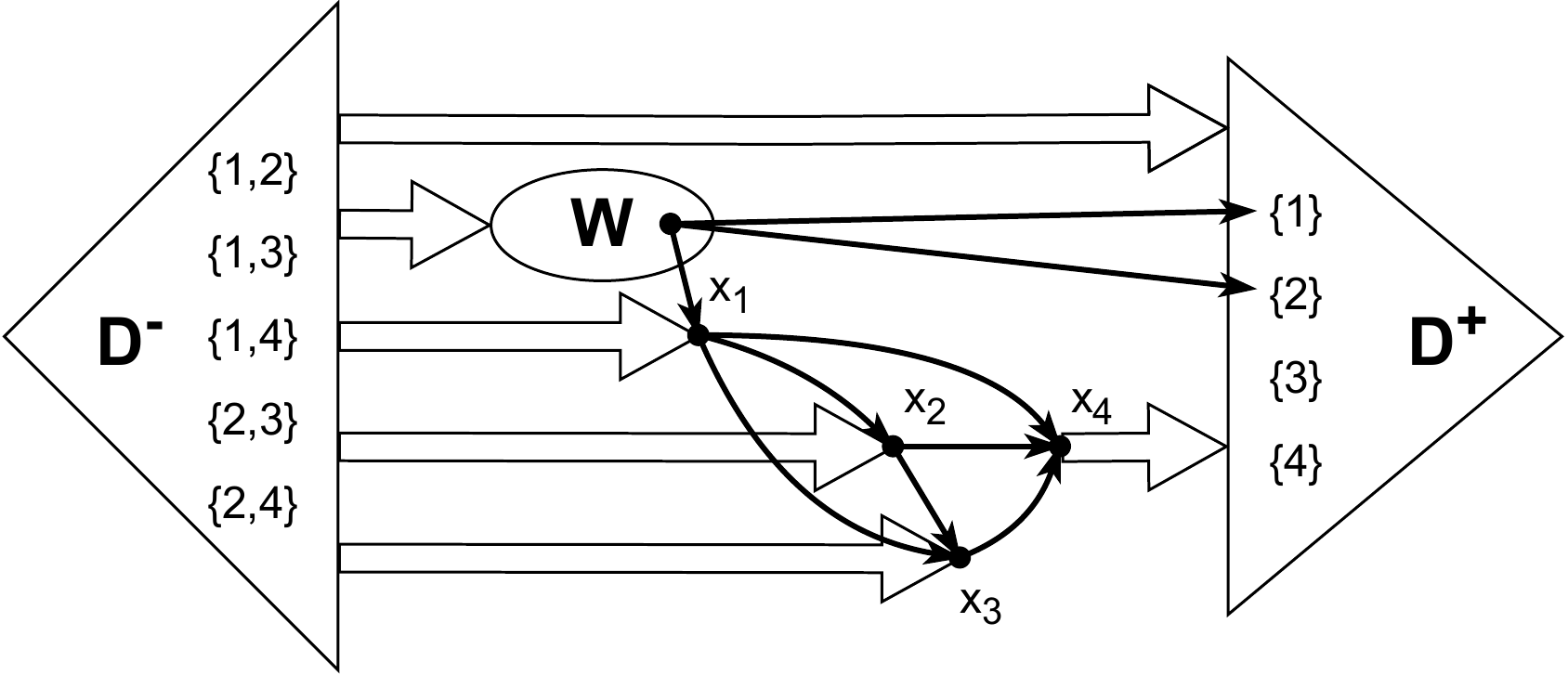}
	\caption{The digraph $D\in\mathcal{D}(4,3)$ that is not covered by $4$ cuts. The tall arrows indicate that there are all possible arcs from the left set to the right set. We find specific characteristic sets in the subdigraphs $D^-$ and $D^+$.}
	\label{c43}
\end{center}
\end{figure} \vspace{-8mm}

Assume that $D$ is covered by $4$ cuts. Then, the characteristic sets of the vertices in $V^-$ and $V^+$ each contain a family of five or four pairwise incomparable subsets of $\{1,2,3,4\}$, respectively, and the two families are disjoint. A family of five pairwise incomparable subsets cannot contain a set of size lower or greater than 2. Without loss of generality, we assume that $V^-$ contains vertices with characteristic sets $\{1,2\},\{1,3\},\{2,3\}, \{1,4\}, \{2,4\}$. It follows that $V^+$ contains vertices with characteristic sets $\{1\},\{2\},\{3\},\{4\}$. Then, $C(x_4)$ must contain at least two elements. Hence, among the three sets $C(x_1),C(x_2),C(x_3) \not\subseteq C(x_4)$ are at most two of size 1. 

Each of the sets $C(x_1),C(x_2),C(x_3)$ and $C(w),\ w\in W,$ contains one element or is equal to $\{3,4\}$. So $C(x)=\{3,4\}$, for some $x\in \{x_1,x_2,x_3\}$. Let $y,z \in V^+$ with $C(y)=\{1\}$ and $C(z)=\{2\}$. The characteristic set of the vertex $w_{\{x,y,z\}} \in W$ contains at least two elements and is not equal to $C(x)=\{3,4\}$, and this is a contradiction. \hfill $\Box$\\

Unfortunately, in the proof of $c(4,3)>4$, we cannot replace $D^-$ by the $4$-th power of some directed path, since this digraph requires only three characteristic sets of size $2$: For example, use the characteristic sets $\dots,\{1,2,4\},\{1,2\},\{1,3,4\},\{1,3\},\{2,3,4\},\{2,3\},\dots$ cyclically. 
We only could reduce the number of layers of $D^-$, but this construction still demands vertices of large total degree (sum of indegree and outdegree). 

One can find a digraph with maximum total degree $10$, and hence in $\mathcal{D}(5,4)$, that is not covered by $4$ cuts. On the other hand, every digraph with maximum total degree at most $8$ is contained in $\mathcal{D}(5,2)$, so it is covered by $4$ cuts.

The proofs of the upper bounds can all be generalized, but this only yields $c(\Delta^-,\Delta^+) \leq k$ for $\Delta^-+\Delta^+=(1+o(1))M(k)$, which is not substantially better than the trivial bound.

With ideas similar to the proof of $c(4,3)>4$, it should be possible to prove $c(7,7)>5$. But for large $k$, it is not clear whether there is a generalization at all. A main barrier is created by the two families $\{B \subseteq \{1,\dots,k\}: 1 \in B, |B|\leq \lfloor k/2 \rfloor\}$ and $\{B \subseteq \{1,\dots,k\}: 1 \notin B, |B|\geq \lfloor k/2 \rfloor\}$ of size $\sim\frac{1}{4} 2^k$, whose union contains $M(k)$ pairwise incomparable sets, and moreover, each set of the first family is incomparable with each set of the second family. \\ \newpage

\section{Complexity of the cut cover problem in acyclic digraphs} \label{Complexity}

In this Section, we study the algorithmic complexity of the decision problem whether an acyclic digraph is covered by $k$ cuts. For $k\leq 2$, the decision can be done efficiently (for $k=2$ by Theorem \ref{or}.\ref{2-cut cover}). For $k\geq 3$, the problem turns out to be NP-complete, even under some very specific additional conditions.

For a digraph $D=(V,A)$, we consider the class $\mathcal{H}(D)$ of all digraphs $D'=(V',A')$ that admit a homomorphism to $D$, that is, a labeling $l:V'\to V$ such that $l(u)l(v) \in A$ for every arc $uv \in A'$. If $D$ is covered by $k$ cuts, then the condition $D'\in \mathcal{H}(D)$ is sufficient for a $k$-cut cover of $D'$. If $D$ is acyclic, then every $D'\in \mathcal{H}(D)$ is acyclic.

We first show that the statement $D \in \mathcal{H}(P_n^d)$ can be decided in polynomial time. $D=(V,A)$ is contained in $\mathcal{H}(P_n^d)$ if and only if there is a labeling $l:V\to\{0,\dots,n-1\}$ such that $l(v) - l(u) \in \{1,\dots,d\}$ for every arc $uv\in A$. We consider only antisymmetric digraphs, since this condition is necessary for being acyclic. An \textit{undirected walk} $U$ in an antisymmetric digraph $D=(V,A)$ is a sequence $(u_0, u_1, \dots, u_m)$ of vertices such that, for every $i\in \{1,\dots,m\}$, there is either a \textit{forward} arc $u_{i-1}u_i$ or a \textit{backward} arc $u_i u_{i-1}$. Let $m^+$ and $m^-=m-m^+$ be the number of forward arcs and backward arcs in $U$, respectively. We define $\Delta(U) := m^+ - dm^-$ depending on $d$.

By proving the following theorem, we extend some ideas that have already been introduced by Brewster and Hell in \cite{homomorphic}.

\textbf{Theorem \ref{Complexity}.\newTh{Hom}\ -\ Testing the existence of a homomorphism to \boldmath $P_n^d$\\ } 
Let $n,d \in \mathbb{N}$. A digraph $D=(V,A)$ is contained in $\mathcal{H}(P_n^d)$ if and only if $\Delta(U) \leq n-1$ for every undirected walk $U$ in $D$. In this case, a homomorphism $l$ from $D$ to $P_n^d$ is defined by
$$l(v) = \max \{ \Delta(U): U \text{ undirected walk in } D \text{ with last vertex } v \}\ ,$$
and can be computed in $O(|V|\cdot |A|)$.

\textbf{Proof. \ } Assume first that $D \in \mathcal{H}(P_n^d)$. Then, there exists a labeling $l:V\to\{0,\dots,n-1\}$ such that $l(v) - l(u) \in \{1,\dots,d\}$ for every arc $uv\in A$. Let $U=(u_0, u_1, \dots, u_m)$ be an undirected walk in $D$ and let $m^+$ and $m^-$ be the number of forward arcs and backward arcs arcs in $U$, respectively. Then,\\ %\vspace*{-3mm}
\hspace*{11mm} $ \begin{aligned} 
 n-1 \ & \geq \ l(u_m)-l(u_0)\ =\ \sum_{i=1}^m \, l(u_i) - l(u_{i-1})\\
\ & \geq \ \hspace*{-10mm} \sum_{\hspace*{9mm} u_{i-1} u_i \in A}\hspace*{-9mm} l(u_i) - l(u_{i-1})\ - \hspace*{-10mm} \sum_{\hspace*{9mm} u_i u_{i-1} \in A} \hspace*{-9mm} l(u_{i-1}) - l(u_i)\ \geq \ m^+ \cdot 1 - m^- \cdot d \ =\ \Delta(U)\,. \end{aligned}$

Now assume that $\Delta(U) \leq n-1$ for every undirected walk $U$ in $D$. We have to prove $D \in \mathcal{H}(P_n^d)$.
Define $l:V\to \{0,\dots,n-1\}$ as given in the statement of the theorem. Indeed for all $v\in V$, $l(v) \geq 0$ by considering the walk $U=(v)$, and $l(v) \leq n-1$ by the assumption.

Let $vw \in A$. Since every undirected walk with last vertex $v$ can be expanded to an undirected walk with last vertex $w$ by adding the forward arc $vw$, it holds that $l(w) \geq l(v) +1$. Similarly, since every undirected walk with last vertex $w$ can be expanded to an undirected walk with last vertex $v$ by adding the backward arc $vw$, it holds that $l(v) \geq l(w) -d$. Combining both statements, we have $l(w) - l(v) \in \{1,\dots,d\}$. Hence, $l$ is a homomorphism to $P_n^d$.\\

Finally, we have to prove that $l$ can be computed in $O(|V|\cdot |A|)$. We can assume that $D=(V,A)$ is weakly connected. For any $v\in V$, let $l^1(v)=0$ and, for $m\geq 2$,

\hspace*{11mm} $l^{m}(v) = \max\left(\ \{0\} \cup \{l^{m-1}(u)+1: uv \in A\} \cup \{l^{m-1}(w)-d: vw \in A\}\ \right) .$

We prove inductively that, for all $v \in V$,

\hspace*{11mm} $l^m(v) = \max \{ \Delta(U): U$ undir. walk in $D$ on at most $m$ vertices with last vertex $v \}.$ 

This is true for $m=1$. Assuming the statement for $m-1$, we conclude for $m\geq 2$

\hspace*{11mm}  $\begin{aligned} l^{m}(v) & = \max\left(\ \{0\} \cup \{l^{m-1}(u)+1: uv \in A\} \cup \{l^{m-1}(w)-d: vw \in A\}\ \right)\\  
\ & = \max\left(\ \{\,\Delta(\,(v)\,)\,\}^{\textcolor{white}{1}} \cup \right. \\
\ & \quad \ \{\Delta(U): U \text{ on at most } m \text{ vertices, forward last arc and last vertex } v\}\ \cup\\
\ & \quad \left.\{\Delta(U): U \text{ on at most } m \text{ vertices, backward last arc and last vertex } v\}^{\textcolor{white}{1}}\right)\\
\ & = \max\{\Delta(U): U \text{ on at most } m \text{ vertices and last vertex } v\}. \end{aligned}$\\

We prove that we only have to determine $l^1,\dots,l^{|V|+1}$. In every iteration, each arc is considered exactly twice. Hence, an iteration can be implemented in $O(|A|)$ (since $D$ is weakly connected, $|V| \leq |A|+1$). Now, two distinct cases could happen: If $l^{|V|} = l^{|V|+1}$, by using the iteration equality inductively, $l^{|V|} = l^m$ for all $m\geq |V|$. In this case, $l = \max_{m=1}^\infty l^m = l^{|V|}$. (We could stop the iteration already at the point, where $l^m$ stabilizes.)

Otherwise, there is a vertex $v$ with $l^{|V|}(v) < l^{|V|+1}(v)$. Then, there exists an undirected walk $U=(v_0,\dots,v_{|V|})$ with last vertex $v_{|V|}=v$ such that $\Delta(U) > \Delta(U')$ for every undirected walk $U'$ on at most $|V|$ vertices with last vertex $v$. 

Two of the $|V|+1$ vertices of $U$ must coincide, say ${u_s=u_t}$ for $0 \leq s < t \leq |V|$. Let $\tilde{U}=(u_s, u_{s+1}, \dots, u_t)$, and let $U'=(u_0,\dots,u_{s-1},u_s,u_{t+1},\dots,u_{|V|})$ be the undirected walk that is obtained from $U$ by skipping $\tilde{U}$. By the assumption, $\Delta(U') < \Delta(U) = \Delta(U') + \Delta(\tilde{U})$. It follows that $\Delta(\tilde{U}) > 0$. For every vertex $u$, we can build undirected walks $U$ with last vertex $u$ and arbitrarily large $\Delta(U)$: We first walk $\tilde{U}$ arbitrarily often and then walk to $u$, using that $D$ is weakly connected. This proves $l(u) = \infty$ for all $u\in V$.\hfill $\Box$\\

\textbf{Theorem \ref{Complexity}.\newTh{NPC Hom}\ -\ NP-Completeness of the cut cover problem in \boldmath $\mathcal{H}(D)$ for acyclic \nolinebreak $D$}\\ 
Let $k\geq 3$, and let $D$ be an acyclic digraph that is not covered by $k$ cuts. Then, the decision problem whether a digraph $D'\in \mathcal{H}(D)$ is covered by $k$ cuts is NP-complete.

\textbf{Proof. \ } Let $D=(W,A)$. Among all digraphs in $\mathcal{H}(D)$ that are not covered by $k$ cuts, let $\tilde{D}$ be one with a minimum number of pairs $(xy,yz)$ of consecutive arcs. We can replace $D$ by $\tilde{D}$ in the claim since all digraphs in $\mathcal{H}(\tilde{D})$ are also in $\mathcal{H}(D)$. 

This has the following effect: We can assume that all $\tilde{D} \in \mathcal{H}(D)$ that have a lower number of pairs $(xy,yz)$ of consecutive arcs than $D$ are covered by $k$ cuts. We now choose a specific $\tilde{D}$. Consider a longest path $\dots,x,y,z$ in $D$ and consider its last three vertices. Then, $y$ has at least one out-neighbor and all of its out-neighbors have outdegree $0$. Let $y'\notin W$, and let $\tilde{D} = (W\cup\{y'\},\{xy'\}\cup A\setminus\{xy\})$ be the digraph that results from $D$ by cutting off $xy$ from $y$, that is, adding a new vertex $y'$ and replacing $xy$ by $xy'$. We prove the following two properties.
\begin{enumerate}
	\item[(1)] For every $k$-cut cover of $\tilde{D}$, $|C(x)| = 1$ and $C(x)\subseteq C(y)$. \vspace*{-2mm}
	\item[(2)] For every $a\in \{1,\dots,k\}$, there is a $k$-cut cover of $\tilde{D}$ with $C(x)=C(y)=\{a\}$.
\end{enumerate}
(1): \ We can assume that the characteristic sets of all vertices $z$ with outdegree $0$ have characteristic set $C(z)=\emptyset$. $C(y)\neq \emptyset$, and since now $y$ has only out-neighbors with empty characteristic sets, we can reduce $C(y)$ to a set $C'(y)\subseteq C(y)$ of size $1$ and still maintain a $k$-cut cover of $\tilde{D}$.
When we add the arc $xy$, we obtain a superdigraph of $D$, which is not covered by $k$ cuts, so $C(x) \subseteq C'(y)$. This implies $C(x) \subseteq C(y)$. Moreover, $|C(x)| \leq 1$, so, since there is the arc $xy'$, $|C(x)| = 1$.

(2): \ Clearly $\tilde{D} \in \mathcal{H}(D)$ and $\tilde{D}$ has a lower number of pairs of consecutive arcs than $D$: Every pair $(wx,xy')$ in $\tilde{D}$ corresponds to a pair $(wx,xy)$ in $D$, but the pair $(xy,yz)$ gets eliminated. Hence, by the assumption, $\tilde{D}$ is covered by $k$ cuts. By (1), $|C(x)|=1$, so by using a permutation of the numbers $\{1,\dots,k\}$, we can achieve $C(x)=\{a\}$ for an arbitrary $a\in\{1,\dots,k\}$. Similarly to the proof of (1), we can reduce $C(y)$ to its subset $C(x)=\{a\}$. This proves (2).\newpage

We now reduce the $k$-colorability problem for graphs to our problem. Let $G=(V,E)$ be a graph. We have to give a construction in polynomial time of an acyclic digraph $D' \in \mathcal{H}(D)$, such that $D'$ is covered by $k$ cuts if and only if $G$ is $k$-colorable.

For every vertex $v\in V$, we take a copy of $\tilde{D}$, that is, a digraph $D_v$ on vertex set $\{w_v=(w,v): w\in W\cup\{y'\}\}$ and arc set $\lbrace w_vz_v: wz\in \{xy'\}\cup A \setminus \{xy\}\rbrace$. Moreover, for every edge $uv\in E$, we add one of the arcs $x_uy_v$ and $x_vy_u$. Let $D'$ be the resulting digraph. Clearly, $D'\in \mathcal{H}(D)$, which is verified by the homomorphism $w_v \mapsto w, y'_v\mapsto y$.

Assume first that $G$ has a $k$-coloring $f:V\to \{1,\dots,k\}$. By (2), for every $v\in V$, there is a $k$-cut cover of $D_v$ with $C(x_v)=C(y_v)=\{f(v)\}$. This yields a $k$-cut cover of the whole digraph $D$: The only arcs we have to check are of the kind $x_uy_v$ for any edge $uv \in E$. By the coloring property, $C(x_u)=\{f(u)\}\not\subseteq \{f(v)\}=C(y_v)$.

Assume now that $D'$ is covered by $k$ cuts. Consider a vertex $u\in V$. By (1), we have $|C(x_u)| =1$, so we can choose the color $f(u)\in\{1,\dots,k\}$ as the unique element of $C(x_u)$. We have to check that $f$ is indeed a coloring of $D'$. Let $uv\in E$, so without loss of generality, we find the arc $x_uy_v$. Again by (1), $C(x_v)\subseteq C(y_v)$. So $C(x_u)\not\subseteq C(y_v)$ implies that $\{f(u)\}=C(x_u)\not\subseteq C(x_v)=\{f(v)\}$, and hence $f(u)\neq f(v)$.  \hfill $\Box$\\ 

In particular, for $k=3$, we obtain the NP-completeness in $\mathcal{H}(P_9^4)$ and $\mathcal{H}(P_{11}^3)$ with Theorem \ref{indegree}.\ref{delta}(iii) and (iv). We can even add some more restrictions. Note that every planar digraph is covered by $4$ cuts since it is (greedily) $M(4)$-colorable.

\textbf{Theorem \ref{Complexity}.\newTh{NP3}\ -\ NP-completeness of the 3-cut cover problem in acyclic and planar digraphs with maximum indegree and outdegree 3\\}
The decision problem whether a digraph is covered by $3$ cuts is NP-complete, even for acyclic and planar digraphs with maximum indegree and outdegree $3$ that are
\begin{enumerate}
\item[(i)] in $\mathcal{H}(P_{12}^3)$, or
\item[(ii)] in $\mathcal{H}(P_{14}^3)$ and in $\mathcal{D}(2,2)$, that is, for every vertex, the sum of its indegree and outdegree is at most 5.
\end{enumerate}
\textbf{Proof. \ } We first study specific digraphs which we are going to use in the construction. We refer to vertices only by their label, which is a natural number. Later, the labels are used as the homomorphism to $P^3_{12}$ or $P^3_{14}$, respectively. In this prove, the labelings are inverted, that is, if an arc goes from a vertex with label $j$ to a vertex with label $i$, then $j-i\in\{1,2,3\}$.

\begin{figure}[H]
\begin{center}
	%\vspace*{-3mm}
	\includegraphics[width=0.9\textwidth]{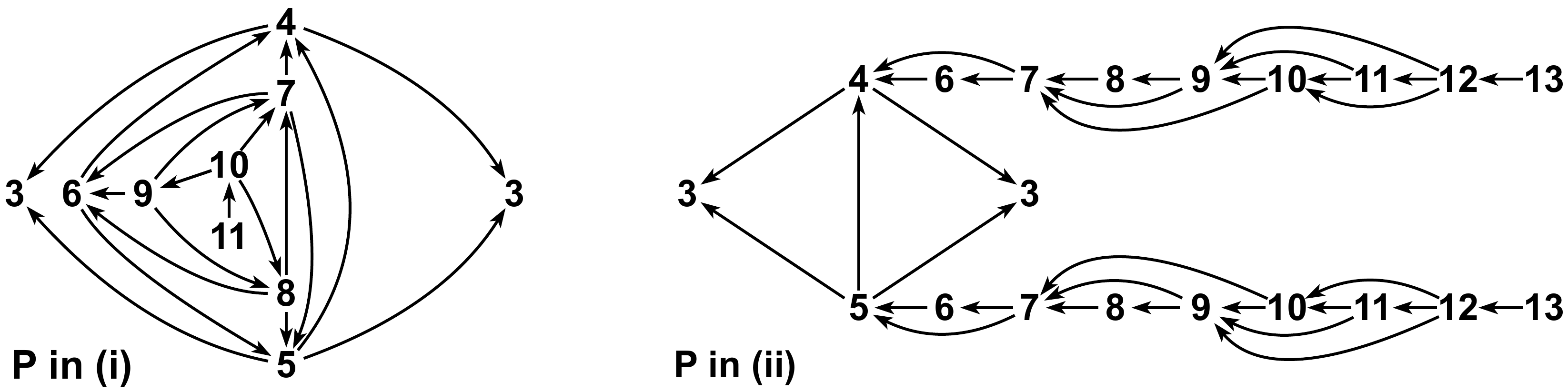}
	\caption{Planar embeddings of the two versions of $P$.}
	\label{P}
\end{center}
\end{figure}

First, we consider a digraph $P$. The different constructions for (i) and (ii) are shown in Figure \ref{P}. We prove the following properties
\begin{enumerate}
\item[(P1)] For every $a \in \{1,2,3\}$, $P$ has a $3$-cut cover with $C(x)=\{a\}$ for both vertices $x$ with \linebreak label $3$.
\item[(P2)] For every 3-cut cover of $P$ with $|C(x)|\geq 1$ for both vertices $x$ with label $3$, there is an $a\in \{1,2,3\}$ such that $C(x)=\{a\}$ for both vertices $x$ with label $3$.
\end{enumerate}

 \vspace{-9mm}

(i): \ For (P1), we use a permutation $\{a,b,c\}=\{1,2,3\}$, and we choose the following sequence of characteristic sets $C(3),C(4),\dots,C(11)$ depending only on the labels of the vertices:
$$\{a\}, \{b\}, \{c\}, \{b,c\}, \{a\}, \{a,b\}, \{a,c\}, \{b,c\}, \{a,b,c\}$$
To prove (P2), consider a 3-cut cover of $P$. It suffices to show that both $C(4)\neq C(5)$ contain exactly one element. Assume that $|C(5)|\geq 2$, say $C(5)=\{2,3\}$. The sets $C(6),C(7),C(8)$ each contain one or two elements and are no subsets of $C(5)=\{2,3\}$. The only such sets are $\{1\}, \{1,2\}, \{1,3\}$. Hence, $C(7)$ and $C(8)$ are the sets $\{1,2\}$ and  $\{1,3\}$. It follows that $C(9)=\{2,3\}$. But then, there is no set of two elements left for $C(10)$. A similar contradictions occurs if $|C(4)|\geq 2$.

(ii): \ For (P1), we use a permutation $\{a,b,c\}=\{1,2,3\}$, and we choose the following sequence of characteristic sets $C(3),C(4),\dots,C(13)$ depending only on the labels of the vertices:
$$\{a\}, \{b\}, \{c\}, \{b,c\}, \{a\}, \{b\}, \{c\}, \{a,b\}, \{a,c\}, \{b,c\}, \{a,b,c\}$$
To prove (P2), consider a 3-cut cover of $P$. It suffices to show that both $C(4)\neq C(5)$ contain exactly one element. Assume that $|C(5)|\geq 2$, say $C(5)=\{1,2\}$. $C(6)$ and $C(7)$ are no subsets of $C(5)$ and contain at most two elements. The only such sets are $\{3\},\{1,3\},\{2,3\}$. Hence, $C(7)$ must also contain two elements, say $C(7)=\{2,3\}$. The characteristic sets $C(8),C(9),C(10)$ contain at most two elements and are no subsets of $C(7)=\{2,3\}$. They can be the sets $\{1\},\{1,2\},\{1,3\}$. It follows that $C(9),C(10)$ are the sets $\{1,2\},\{1,3\}$. Hence, $C(11)=\{2,3\}$. But then, there is no set of two elements left for $C(12)$, and this is a contradiction. So $|C(5)|=1$, and similarly $|C(4)|=1$.

\begin{figure}[H]
\begin{center}
	%\vspace*{-3mm}
	\includegraphics[width=40mm]{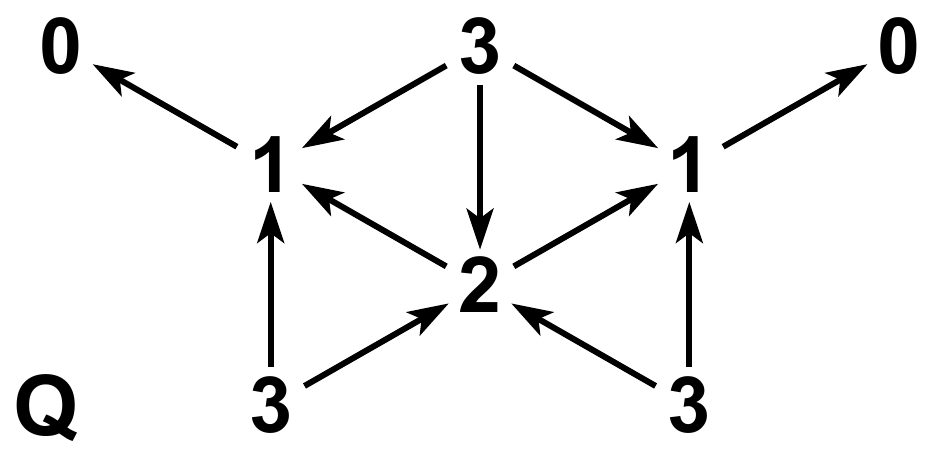}
	\caption{Planar embedding of the digraph $Q$.}
	\label{Q}
\end{center}
\end{figure}

Let $Q$ be the digraph shown in Figure \ref{Q}. $Q$ satisfies the following properties
\begin{enumerate}
\item[(Q1)] For every $a \in \{1,2,3\}$, $Q$ has a $3$-cut cover with $C(x)=\{a\}$ for all vertices $x$ with \linebreak label $3$.
\item[(Q2)] For every 3-cut cover of $Q$ with $|C(x)|=1$ for all vertices $x$ with label $3$, there is an $a\in \{1,2,3\}$ such that $C(x)=\{a\}$ for all vertices $x$ with label $3$.
\end{enumerate}

 \vspace{-8mm}
For (Q1), we use a permutation $\{a,b,c\}=\{1,2,3\}$, and we choose the characteristic sets depending on the label of the vertices: $\emptyset$ for label 0, $\{c\}$ for label 1, $\{b\}$ for label 2, and $\{a\}$ for label 3.

For (Q2), assume first that the characteristic set of the vertex with label 2 contains two elements. Then, the characteristic sets of the vertices with label 3 contain exactly the third element of $\{1,2,3\}$, so all of them are equal. Otherwise, the vertices with label 1 have two ingoing arcs from vertices with distinct characteristic sets of one element. Hence, their characteristic sets also contain one element. So the characteristic sets are a 3-coloring of the vertices with label 1, 2, 3. It follows that the characteristic sets of the vertices with label 3 must coincide.\\

We now reduce the $3$-colorability problem for planar graphs to our problem. Let ${G=(V,E)}$ be a planar graph. We have to give a polynomial construction of a digraph ${D=(W,A)}$ with the required properties that is covered by $3$ cuts if and only if $G$ is $3$-colorable. All vertices in $W$ get labels in $\{0,\dots,11\}$ or $\{0,\dots,13\}$ that prove $D \in \mathcal{H}(P_{12}^3)$ or $D \in \mathcal{H}(P_{14}^3)$, respectively. To avoid a complex naming, we refer to vertices only with their label.

Let $v\in V$ be a vertex of $G$, and let $d(v)\geq 3$ be the degree of $v$. We describe the construction of the subdigraph $D_v$ of $D$ that corresponds to $v$. Let $T_v$ be an arbitrary tree with $d(v)$ vertices of degree $1$ and $d(v)-2$ vertices of degree $3$. We replace the vertices of $T_v$ of degree $3$ by copies of $Q$ and the edges of $T_v$ by copies of $P$. We glue them together in the vertices with label $3$. Let $W_v$ be the set of the $d(v)$ vertices with label $3$ that are not contained in a copy of $Q$. For each of them, we add a vertex with label $0$ and an arc to this vertex. This yields the subdigraph $D_v$ of $D$ that corresponds to $v$. \newpage
\begin{figure}[H]
\begin{center}
	\vspace*{-4mm}
	\includegraphics[height=9cm]{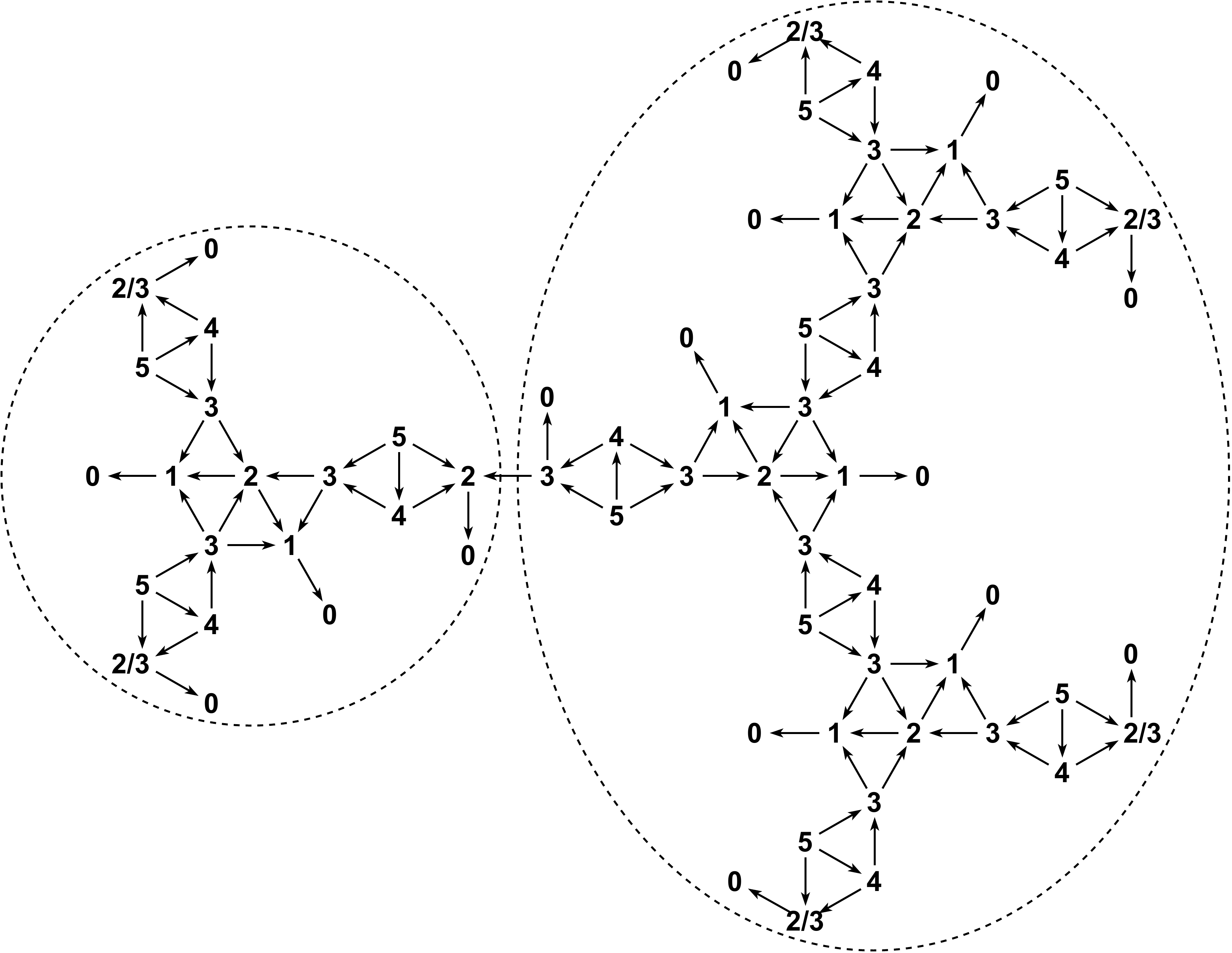}
	\caption{The subdigraph of $D$ that corresponds to two adjacent vertices in $G$, the left one of degree $3$, and the right one of degree $5$. The vertices with label greater than 5 are not shown.}
	\label{NPC}
\end{center}
\end{figure} \vspace{-8mm}
Assume that $D_v$ is covered by $3$ cuts. Since all vertices in $D_v$ with label $3$ have an outgoing arc, their characteristic sets are non-empty. So by (P2), they contain exactly one element, and by (P2) and (Q2), all of their characteristic sets must be equal. By (P1) and (Q1), there is also one such 3-cut cover. This yields the following
\begin{itemize}
\item[(W1)] For every $a\in\{1,2,3\}$, $D_v$ has a $3$-cut cover with $C(w)=\{a\}$ for all $w\in W_v$.
\item[(W2)] For every $3$-cut cover of $D_v$, there is an $a\in\{1,2,3\}$ such that $C(w)=\{a\}$ for all $w\in W_v$.
\end{itemize}

For each edge $uv\in E$, we choose a vertex $w_{uv}\in W_u$ and a vertex $w_{vu} \in W_v$ such that no vertex is chosen twice and that we can draw an arc between $w_{uv}$ and $w_{vu}$ for all $uv \in E$ in a planar way. It is intuitively true that a planar construction is possible, and we do not give a formal proof. To maintain the property of our labeling, we have to relabel one of the vertices $w_{uv},w_{vu}$ from $3$ to $2$. Let $D$ be the resulting digraph. 

Clearly, the labeling of the vertices verifies $D\in \mathcal{H}(12,3)$ or $D\in \mathcal{H}(14,3)$, respectively. In particular, $D$ is acyclic. Moreover, every vertex of $D$ has outdegree and indegree at most 3, and in (ii), not both are equal to 3. Lastly, we have to prove that $G$ is 3-colorable if and only if $D$ is covered by $3$ cuts.

Assume first that $G$ has a 3-coloring $f:V\to\{1,2,3\}$. For every $v\in V$, by (W1), there exists a $3$-cut cover of $D_v$ with $C(w)=\{f(v)\}$ for all $w\in W_v$. These characteristic sets define a 3-cut cover of $D$. The only arcs we have to check are of the form $w_{uv}w_{vu}$ for some edge $uv\in E$. By the coloring property, $C(w_{uv}) = \{f(u)\} \not\subseteq \{f(v)\} = C(w_{vu})$.

Assume now that $D$ is covered by $3$ cuts. Let $v\in V$. By (W1), there is an $a\in\{1,2,3\}$ such that $C(w)=\{a\}$ for all $w\in W_v$. We choose $f(v)=a$. The labeling $f:V\to\{1,2,3\}$ is a coloring of $G$: For every edge $uv\in E$, there is an arc between the vertices $w_{uv}$ and $w_{vu}$ and this yields $\{f(u)\} = C(w_{uv}) \neq C(w_{vu}) = \{f(v)\}$. Hence $f(u)\neq f(v)$. \hfill $\Box$  \newpage

\section{Appendix} \label{Appendix}

The following two bounds on $M(k)$ can be proved simply by induction. 

\textbf{Theorem \ref{Appendix}.\newTh{approxM}\ -\ Simple bounds on \boldmath $M(k)$ } 
$$\frac{\sqrt{2}}{2} \frac{2^k}{\sqrt{k+1}} \ \leq\  M(k)\ \leq\ \frac{\sqrt{3}}{2} \frac{2^k}{\sqrt{k+1}}\ .$$
\textbf{Proof. \ } By induction, we prove the following stronger bounds.

If $k$ is odd, \ $\frac{\sqrt{2}}{2} \frac{2^k}{\sqrt{k+1}} \ \leq\  M(k)\ \leq\ \frac{\sqrt{3}}{2} \frac{2^k}{\sqrt{k+2}}$. \
If $k$ is even, \ $\frac{\sqrt{2}}{2} \frac{2^k}{\sqrt{k}} \ \leq\  M(k)\ \leq\ \frac{\sqrt{3}}{2} \frac{2^k}{\sqrt{k+1}}$.

The bounds are true for $k=1$. So let $k\geq 2$, and assume the bounds for $k-1$.

If $k=2l$ is even, we have that $M(k) = \binom{2l}{l} = \frac{(2l)!}{l!\ l!} = \frac{2l}{l}\frac{(2l-1)!}{(l-1)!\ l!} = 2\binom{2l-1}{l-1} = 2M(k-1)$, so the claim follows directly from the induction hypothesis.
So let $k=2l+1$ be odd. We have that $M(k) = \binom{2l+1}{l} = \frac{(2l+1)!}{(l+1)!\ l!} = \frac{2l+1}{l+1}\frac{(2l)!}{l!\ l!} = \frac{2l+1}{2l+2}2\binom{2l}{l} = \frac{k}{k+1}2\,M(k-1)$.

The induction hypothesis yields for the upper bound,\\
\hspace*{1cm} $M(k) \leq \frac{k}{k+1} 2 \frac{\sqrt{3}}{2} \frac{2^{k-1}}{\sqrt{k}} = \sqrt{\frac{k(k+2)}{(k+1)^2}} \frac{\sqrt{3}}{2} \frac{2^{k}}{\sqrt{k+2}} < \frac{\sqrt{3}}{2} \frac{2^{k}}{\sqrt{k+2}}$,

and similarly, for the lower bound,\\
\hspace*{1cm} $M(k) \geq \frac{k}{k+1} 2 \frac{\sqrt{2}}{2} \frac{2^{k-1}}{\sqrt{k-1}} = \sqrt{\frac{k^2}{(k+1)(k-1)}} \frac{\sqrt{2}}{2} \frac{2^k}{\sqrt{k+1}} > \frac{\sqrt{2}}{2} \frac{2^k}{\sqrt{k+1}}$.  \hfill $\Box$\vspace*{3mm}

Finally, some thoughts on cut covers in symmetric digraphs. An alternative proof for the following theorem of Poljak and Rödl \cite{symmetric} is presented.

\textbf{Theorem \ref{Appendix}.\newTh{symmetric}\ -\ Cut cover of symmetric digraphs \cite{symmetric}}\\
A symmetric digraph $D$ is covered by $k$ cuts if and only if it is $M(k)$-colorable.

\textbf{Proof. \ } We only need to prove one direction. Assume that $D$ is covered by $k$ cuts and consider its characteristic sets. We follow the ideas of a proof of Sperner's theorem that is presented in \cite{Sperner}. (In \cite{symmetric}, another proof of the present theorem is given.)

Choose the $k$-cut cover such that $\sum_{v \in V} \left | |C(v)| - \lfloor k/2 \rfloor \right |$ is minimized. If the sum is $0$, all characteristic sets have the size $\lfloor k/2 \rfloor$, and hence, they give a coloring with $M(k)$ colors. Otherwise, there is a set of size either greater or smaller than $\lfloor k/2 \rfloor$. We consider the first case. The second one is done similarly. 

Let $s=\max \{ |C(v)| : \ v\in V \} > \lfloor k/2 \rfloor$.
Our goal is to replace each set $C \in \mathcal{P}_s \linebreak = \{C \subseteq \{1,\dots,k\}: \ |C| = s\}$ by a set $B(C) \in \mathcal{B}_C = \{B \in \mathcal{P}_{s-1}: B \subseteq C\}$ such that $B(C) \neq B(C')$ for each two $C \neq C'$ in $\mathcal{P}_s$. We are able to do this if and only if $(\mathcal{B}_C)_{C \in \mathcal{P}_s}$ has a system of distinct representatives. And this is true by Hall's theorem: Let $\mathcal{C} \subseteq \mathcal{P}_s$ and $\mathcal{B} = \bigcup_{C \in \mathcal{C}} \mathcal{B}_C$. By double-counting and using that $2s \geq k+1$, we obtain 
$$|\mathcal{C}| \cdot s\, =\, |\{ \, (B,C): C\in \mathcal{C}, B\in \mathcal{B}_{C} \}| \ \leq \ |\mathcal{B}|\cdot (k-s+1)\, \leq \, |\mathcal{B}|\cdot s.$$
Dividing by $s$ yields Hall's condition $|\mathcal{C}| \leq |\mathcal{B}|$.

For $C\notin \mathcal{P}_s$, let $B(C)=C$, so in every case $B(C)\subseteq C$. We now replace the characteristic sets $C(v)$ by $B(C(v))$. Let $uv \in A$. We have to show $B(C(u)) \not\subseteq B(C(v))$. So assume that $B(C(u)) \subseteq B(C(v))$. 
In the case $C(u) \notin \mathcal{P}_s$, we have $C(u) = B(C(u)) \subseteq B(C(v)) \subseteq C(v)$, and this contradicts $C(u) \not\subseteq C(v)$.

In the case $C(u) \in \mathcal{P}_s$, we have that $|B(C(u))| = s-1 \geq |B(C(v))|$, so $B(C(u)) \subseteq B(C(v))$ already implies $B(C(u)) = B(C(v))$. This cannot happen if both $C(u),C(v) \in \mathcal{P}_s$, since $C(u) \neq C(v)$. So let $C(v)\notin \mathcal{P}_s$. It follows that $C(v) = B(C(v)) = B(C(u)) \subseteq C(u)$. Since there is also an arc $vu$, this is a contradiction.

The new characteristic sets $B(C(v))$, $v \in V$, contradict the choice of the $k$-cut cover. \hfill $\Box$\\

This theorem is a generalisation of Sperner's theorem, which we obtain by considering complete digraphs, that is, all characteristic sets must be incomparable.

In \cite{NPC}, Watanabe, Sengoku, Tamura and Shinoda proved that the decision problem whether a digraph is covered by $k$ cuts is NP-complete for all $k\geq 3$, and for $k=3$ even for planar and symmetric digraphs. The previous theorem implies a stronger result.

\textbf{Corollary \ref{Appendix}.\newTh{symmetric NPC}\ -\ NP-completeness for symmetric digraphs \cite{symmetric}\\}
For $k \geq 3$, the decision problem whether a symmetric digraph is covered by $k$ cuts is NP-complete. For $k = 3$, this is still true when restricted on planar digraphs.

\textbf{Proof. \ } The $M(k)$-colorability problem of graphs can be reduced to the $k$-cut cover problem of symmetric digraphs by replacing every edge by two opposite arcs. In the case $k=M(k)=3$, we can add the condition that the graphs and digraphs are planar \cite{symmetric}. \hfill $\Box$

\end{document}